\documentclass[leqno,12pt]{amsart}
\usepackage{latexsym}
\usepackage{mathrsfs}
\usepackage{amsmath}
\usepackage{amssymb}
\usepackage[bookmarks]{hyperref}
\usepackage{pstricks,pst-plot}

\newtheorem{theorem}{Theorem}[section]
\newtheorem{lemma}[theorem]{Lemma}
\newtheorem{corollary}[theorem]{Corollary}
\newtheorem{proposition}[theorem]{Proposition}

\theoremstyle{definition}

\newtheorem{definition}[theorem]{Definition}

\newcommand{\C}{\mathbb{C}}

\newcommand{\R}{\mathbb{R}}

\newcommand{\bthm}{\begin{theorem}}
\newcommand{\ethm}{\end{theorem}}
\newcommand{\blem}{\begin{lemma}}
\newcommand{\elem}{\end{lemma}}
\newcommand{\bcor}{\begin{corollary}}
\newcommand{\ecor}{\end{corollary}}
\newcommand{\bprop}{\begin{proposition}}
\newcommand{\eprop}{\end{proposition}}
\newcommand{\bdefn}{\begin{definition}}
\newcommand{\edefn}{\end{definition}}
\newcommand{\bpf}{\begin{proof}}
\newcommand{\epf}{\end{proof}}

\def\vep {\varepsilon}
\def \sm {\setminus}

\def\od{\overline D}

\def\h#1{\widehat {#1}}
\def\hh#1{\widehat {#1} \sm  {#1}}

\def \ab {{[a_j, b_j]}}
\def \openab {(a_j,b_j)}
\def \abn#1{[a_1, b_1] \cup \cdots \cup [a_{#1}, b_{#1}]}
\def \openabn#1{(a_1, b_1) \cup \cdots \cup (a_{#1}, b_{#1})}

\def \Int {{\rm Int}}
\def\I {[0,1]}

\begin{document}
\title{Polynomial Hulls of Arcs and Curves}
\author{Alexander J. Izzo}
\address{Department of Mathematics and Statistics, Bowling Green State University, Bowling Green, OH 43403}
\email{aizzo@bgsu.edu}
\thanks{The author was supported by NSF Grant DMS-1856010.}


\begin{abstract}
It is shown that there exist arcs and simple closed curves in $\C^3$ with nontrivial polynomial hulls that contain no analytic discs.  It is also shown that in any bounded, connected Runge domain of holomorphy in $\C^N$ ($N \geq 2$) there exist polynomially convex arcs and simple closed curves of almost full measure.  These results, which strengthen earlier results of the author, are obtained as consequences of a general result about polynomial hulls of arcs and simple closed curves through Cantor sets.
\end{abstract}

\maketitle

\section{Introduction}

The main purpose of this paper is to prove that (i) there exist arcs and simple closed curves in $\C^3$ with nontrivial polynomial hulls that contain no analytic discs, and (ii) there exist polynomially convex arcs and simple closed curves in $\C^N$ ($N\geq 2$) that are large in the sense of 
$2N$-dimensional Lebesgue measure.  The precise results, which we state below after a few preliminary remarks on our terminology and notation, give more detailed information.

By an {\it arc} we mean a space homeomorphic to the closed unit interval and by a {\it simple closed curve}, a space homeomorphic to the unit circle.  Thus the arcs and simple closed curves are the compact, connected, one-dimensional manifolds.
Throughout the paper, $m$ will denote $2N$-dimensional Lebesgue measure on $\C^N$.
For a compact set $X$ in $\C^N$, we denote by $P(X)$ the uniform closure on $X$ of the polynomials in the complex coordinate functions $z_1,\ldots, z_N$.
The \emph{polynomial hull} 
$\h X$ of $X$ is defined by
$$\h X=\{z\in\C^N:|p(z)|\leq \max_{x\in X}|p(x)|\ 
\mbox{\rm{for\ all\ polynomials}}\ p\}.$$
The set $X$ is said to be \emph{polynomially convex} if $\h X = X$.  The polynomial hull of $X$ is said to be \emph{nontrivial} if instead the set $\h X \sm X$ is nonempty.
By an \emph{analytic disc} in $\C^N$, we mean an injective holomorphic map $\sigma: \{ z\in \C: |z|< 1\}\rightarrow\C^N$.
By the statement that a subset $S$ of $\C^N$ contains no analytic discs, we mean that there is no analytic disc in $\C^N$ whose image is contained in $S$.

Following Garth Dales and Joel Feinstein \cite{DalesF},
we shall say that a Banach algebra $A$ has dense invertibles if the invertible elements of $A$ are dense in $A$.  As noted in \cite{DalesF}, for $X$ a compact set in $\C^N$, the condition that $P(X)$ has dense invertibles is strictly stronger than the condition that $\h X$ contains no analytic discs.

\bthm\label{arc-no-disc}
There exists an arc $J$ in $\C^3$ such that the polynomial hull $\h J$ of $J$ is strictly larger than $J$ and $P(J)$ has dense invertibles.  The same statement holds with 
\lq\lq arc\rq\rq\ replaced by \lq\lq simple closed curve\rq\rq.
\ethm

\bthm\label{poly-convex}
Let $\Omega$ be a bounded, connected Runge domain of holomorphy in $\C^N$, fix a point $x_0$ in $\Omega$, and fix $\vep>0$.  
Then there exists a polynomially convex arc $J$ in $\C^N$ such that $x_0\in J \subset \Omega$ and $m(\Omega\sm J)<\vep$.  Furthermore, $J$ can be chosen so that $P(J)=C(J)$ and the set of polynomials zero-free on $J$ is dense in 
$P(\overline\Omega)$.  The same statements hold with \lq\lq arc\rq\rq\ replaced by \lq\lq simple closed curve\rq\rq\ provided $N\geq 2$.  
\ethm

It will be shown that these two results follow readily from the existence of certain Cantor sets together with the following general result, to be proven below, about polynomial hulls of arcs and curves through Cantor sets.

\bthm\label{general-theorem}
Let $E$ be a Cantor set in $\C^N$.  Then $E$ is contained in an arc $J$ in $\C^N$ such that 
$\h J= J \cup \h E$.  Furthermore, $J$ can be chosen so that the closure of each component of $J\sm E$ is a $C^\infty$-smooth arc.  If $\Omega$ is a connected Runge domain of holomorphy in $\C^N$ that contains $E$, then $J$ can be chosen to lie in $\Omega$.
The same statements hold with \lq\lq arc\rq\rq\ replaced by \lq\lq simple closed curve\rq\rq\ provided $N\geq 2$.
\ethm

We remark that when $N=1$ the hypothesis that $\Omega$ is Runge will not be used in the proof.

As an almost immediate consequence of Theorem~\ref{general-theorem} we will obtain the following.

\bcor\label{corollary}
Each polynomially convex Cantor set $E$ in $\C^N$ is contained in a polynomially convex arc $J$ in $\C^N$ that satisfies $P(J)=C(J)$.  The same statements hold with \lq\lq arc\rq\rq\ replaced by \lq\lq simple closed curve\rq\rq\ provided $N\geq 2$.
\ecor

A few historical remarks to put Theorems~\ref{arc-no-disc} and~\ref{poly-convex} in context are in order.
Polynomial hulls of arcs and curves have been studied extensively.  John Wermer \cite{Wermer-Cantor}
showed that there exist arcs in $\C^3$ with nontrivial polynomial hull, and 
 by modifying Wermer's construction, Walter Rudin \cite{Rudin} extended the result to arcs in  
$\C^2$.  In contrast, though, Wermer \cite{Wermer1958} also showed that every {\em real-analytic} arc in $\C^N$ is polynomially convex and that the hull of every {\em real-analytic\/} simple closed curve in $\C^N$ is 
a 1-dimensional complex-analytic variety.  The regularity hypothesis was weakened by many mathematicians, and it is now known, by work of Herbert Alexander \cite{A}, that real-analyticity can be replaced by rectifiability.  Theorem~\ref{poly-convex} above shows that some highly {\em nonrectifiable} arcs and simple closed curves are also polynomially convex.
Theorem~\ref{poly-convex} also strengthens \cite[Lemma~3.3]{Izzo2} which gives a {\em rationally convex} connected set with properties similar to those of the arc in Theorem~\ref{poly-convex}.

The existence of analytic structure in polynomial hulls has also been studied extensively.
It was once conjectured that whenever the polynomial hull $\h X$ of a compact set $X$ in $\C^N$ is strictly larger than $X$, the complementary set $\h X\sm X$ must contain an analytic disc. This conjecture was disproved by Gabriel Stolzenberg \cite{Stol1} who gave a counterexample in $\C^2$.  Numerous additional counterexamples have been presented in the literature since then.
The first example of a compact set $X$ in $\C^2$ with nontrivial polynomial hull such that the uniform algebra $P(X)$ has a dense set of invertible elements was given by Dales and Feinstein \cite{DalesF}.
Recent work of the author, H\aa kan Samuelsson Kalm, and Erlend Forn{\ae}ss Wold \cite{ISW} and of the author and Lee Stout \cite{IS} shows that every smooth manifold of dimension strictly greater than one smoothly embeds in some $\C^N$ as a subspace $X$ such that $\h X\sm X$ is nonempty but contains no analytic discs.
(The most recent result in this direction is due to Leandro Arosio and Wold \cite[Corollary~1.3]{AWold}.)
In response to a talk on this work given by the author, Hari Bercovici raised the question of whether a (nonsmooth) {\em one-dimensional} manifold can have a nontrivial polynomial hull containing no analytic discs. In fact a similar question, but specifically for manifolds in $\C^2$, was raised by Wermer \cite{Wermer1954} more than 60 years ago.  The author \cite[Theorems~1.7]{Izzo} answered Bercovici's question affirmatively by giving an example in $\C^4$.  Theorem~\ref{arc-no-disc} strengthens that result by decreasing the dimension of the ambient space.
The question of whether there exists an example in $\C^2$ remains open.  The proof of Theorem~\ref{arc-no-disc} given here shows, however, that if there is a \emph{Cantor set} in $\C^2$ with nontrivial polynomial hull that contains no discs, then there is also an arc in $\C^2$ with this property.

The author would like to thank Lee Stout for very helpful and inspiring correspondence on the topic of this paper.  Indeed much of the content of Theorems~\ref{arc-no-disc} and~\ref{general-theorem} was found by Stout (by a different, but related, argument) and communicated to the author.  Specifically the first assertion of Theorem~\ref{general-theorem} (the existence of an arc $J$ such that $\h J = J \cup \h E$) was found by Stout and used by him to obtain an arc in 
$\C^3$ with nontrivial polynomial hull that contains no analytic discs.
Also consideration of Runge domains in Theorem~\ref{poly-convex}, rather than only the unit ball, was suggested by Stout.

This paper was written while the author was a visitor at the University of Michigan.  He would like to thank the Department of Mathematics for its hospitality.

\section{The proofs}

The results stated in the introduction will be proved in the following order:
Theorem~\ref{general-theorem}, Corollary~\ref{corollary}, Theorem~\ref{poly-convex}, Theorem~\ref{arc-no-disc}.  

We will use the following standard terminology and notation.  Given a number $\alpha>0$ and a set $A$ in $\C^N$, by the {\it $\alpha$-neighborhood\/} of $A$ we mean the set of points in $\C^N$ whose distance from $A$ is strictly less than $\alpha$.  The supremum-norm of a function $f$ will be denoted by $\|f\|_\infty$, and for each $k=1, 2, \ldots$, the $C^k$-norm of $f$ will be denoted by $\|f\|_{C^k}$.
The real part of a complex number (or function) $z$ will be denoted by $\Re z$.

The proof of Theorem~\ref{general-theorem} relies on two results which we quote here for the reader's convenience: a result of Whyburn~\cite{Whyburn} on the existence of an arc through a totally disconnected set, and a result of Arosio and Wold \cite{AWold} concerning perturbing a totally real embedding of a compact manifold so as to control the polynomial hull of its image.

In Whyburn's terminology, a \emph{continuous curve} is a connected, locally connected, locally compact, separable metric space, and a \emph{local separating point} of a continuous curve $M$ is a point that is a cut point of some connected open subset of $M$.

\bthm\cite[Theorem p.~57]{Whyburn}\label{Whyburn-theorem}
If $K$ is any closed, compact and totally disconnected subset of a continuous curve $M$ having no local separating point and $p$ and $q$ are any two points of $K$, then there exists in $M$ an arc $pq$ which contains $K$.
\ethm

\bcor\label{existence}
Let $E$ be a Cantor set contained in a connected open subset $\Omega$ of $\R^N$, $N\geq 2$.  Then there exists an arc in $\Omega$ that contains $E$ and whose end points lie in $E$, and there also exists a simple closed curve in $\Omega$ that contains $E$.
\ecor

\bpf
The existence of the desired arc is obviously just a special case of Whyburn's theorem.  To construct the simple closed curve, choose distinct points $p$ and $q$ in $\Omega\sm E$, and set $K=E\cup \{p,q\}$.  
 By Whyburn's theorem there exists in $\Omega$ an arc $pq$ which contains $K$.  Denote this arc by $J$.
By modifying $J$ near the end points $p$ and $q$, we may assume that there are open Euclidean balls $B_p$ and $B_q$ centered at $p$ and $q$ respectively, and whose closures lie in $\Omega$ such that the intersection of $J$ with each of $B_p$ and $B_q$ is a straight line segment.  Choose points $p'$ and $q'$ in $B_p\sm J$ and $B_q\sm J$ respectively.  The set $\Omega\sm J$ is path connected, so there is a path from $p'$ to $q'$ in $\Omega\sm J$.  By discarding initial and final segments of this path, we can obtain an arc $\widetilde J$ in $\Omega\sm (J\cup B_p \cup B_q)$ whose end points $\tilde p$ and $\tilde q$ lie on the boundary of $B_p$ and $B_q$ respectively.  Let $L_{\tilde p}$ and $L_{\tilde q}$ be the straight line segments $p\tilde p$ and $q\tilde q$ respectively.  Then $J\cup L_{\tilde q} \cup \widetilde J \cup L_{\tilde p}$ is a simple closed curve in $\Omega$ that contains $E$.
\epf

\bthm\cite[Theorem~1.4]{AWold}\label{AWold-theorem}
Let $M$ be a compact $C^\infty$-manifold {\rm(}possibly with boundary{\rm)} of dimension $d<n$, and let $f:M\rightarrow \C^N$ be a totally real $C^\infty$-embedding.  Let $K\subset \C^N$ be a compact polynomially convex set.  Then, for any $k\geq 1$ and for any $\vep>0$, there exists a totally real $C^\infty$-embedding $f_\vep : M\rightarrow \C^N$ such that the following hold:
{\baselineskip=16pt
\begin{enumerate}
\item[(1)] $\| f_\vep - f\|_{C^k} < \vep$
\item[(2)] $f_\vep=f$ on $f^{-1}(K)$
\item[(3)] $f_\vep(M) \cup K$ is polynomially convex.
\end{enumerate}
}
\ethm

Note that every smooth embedding of a one-dimensional manifold is automatically totally real.

We now turn to the proofs of the results stated in the introduction.

\bpf[Proof of Theorem~\ref{general-theorem}]
The case $N=1$ is special and easy; it follows from Corollary~\ref{existence} and the fact that in the plane every Cantor set and every arc is polynomially convex.
(Note that the hypothesis that $\Omega$ is Runge is unneeded when $N=1$.)

From now on we assume that $N\geq 2$.  We treat only the construction of the arc, the construction of the simple closed curve being similar.  
Let $\Omega$ be a connected Runge domain of holomorphy in $\C^N$ that contains $E$.
By Corollary~\ref{existence} there is an arc in $\Omega$
that contains the Cantor set $E$ and whose end points lie in $E$.  
Given the existence of such an arc, 
one can show that there is such an arc $J_0$ with the additional property that the closure of each component of $J_0\sm E$ is a $C^\infty$-smooth arc.  In other words, there is a topological embedding $\sigma_0:\I \rightarrow \C^N$ with $\sigma_0(\I)\subset \Omega$ and $\sigma_0(\{0,1\})\subset E$ such that letting $J_0=\sigma_0(\I)$, letting $K=\sigma_0^{-1}(E)$, and letting $(a_1, b_1), (a_2, b_2), \ldots$ be the components of $\I\sm K$, we have that $J_0$ contains $E=\sigma_0(K)$ and each restriction map $\sigma_0|_\ab \rightarrow \C^N$ is a $C^\infty$-embedding.  

The proof will be complete once we establish that there is a 
topological embedding $\sigma:\I \rightarrow \C^N$ such that $\sigma(\I)\subset \Omega$, such that $\sigma$ agrees with $\sigma_0$ on $K$, such that each restriction map $\sigma|_\ab \rightarrow \C^N$ is a $C^\infty$-embedding, and such that in addition, setting $J=\sigma(\I)$, we have that $\h J = J \cup \h E$.  We will obtain the map $\sigma$ as a uniform limit of a sequence of continuous maps $\sigma_n: \I \rightarrow \C^N$.  The construction of the sequence $(\sigma_n)$ will be by induction.  We will simultaneously also construct a sequence of compact polynomially convex sets $X_0, X_1, \ldots$ and a decreasing sequence of strictly positive numbers $\alpha_0, \alpha_1,\ldots$.

Because $\Omega$ is a Runge domain of holomorphy, $\Omega$ contains the polynomial hull of $J_0\cup \h E$, and hence,
there is a compact polynomially convex neighborhood of $J_0\cup \h E$ contained in $\Omega$, which we will take as $X_0$.  
Choose $\alpha_0>0$ such that $X_0$ contains the $\alpha_0$-neighborhood of $J_0 \cup \h E$.

Before presenting the conditions we will require the sequences $(\sigma_n)$, $(X_n)$, and $(\alpha_n)$ to satisfy, we make some preparations.  First choose, for each $n=1, 2, \ldots$, a compact polynomially convex neighborhood $L_n$ of $\h E$ contained in the $(1/n)$-neighborhood of $\h E$.  Since, for each $j=1,2,\ldots$, the map $\sigma_0|_\ab$ is a $C^\infty$ embedding, the stability of smooth embeddings gives that there exists $\delta_j>0$ such that every $C^\infty$-map $\alpha:\ab \rightarrow \C^N$ satisfying $\bigl\|\alpha - \sigma_0|_\ab\bigr\|_{C^1} <\delta_j$ is a $C^\infty$-embedding.  For $s\in \openab$, let $d(s)$ denote the distance from $\sigma_0(s)$ to $\sigma_0(\I\sm\openab)$.
Then $d$ is a continuous function on $\bigcup_{j=1}^\infty \openab= \I \sm K$ and is everywhere strictly positive.

The sequence of continuous maps $(\sigma_n)$, the sequence of compact polynomially convex sets $(X_n)$, and the decreasing sequence of strictly positive numbers $(\alpha_n)$ will be chosen so that the following conditions hold for all 
$n=1,2,\ldots$, all $j=1,2,\ldots$, and all $s\in \I \sm K$.

{\baselineskip=17pt
\begin{enumerate}
\item[(i)] $\sigma_n|_K=\sigma_0|_K$
\item[(ii)] each map $\sigma_n|_\ab$ is of class $C^\infty$
\item[(iii)]  $\| \sigma_n - \sigma_{n-1}\|_\infty < \alpha_{n-1}/ 2^n$
\item[(iv)] $\bigl | \sigma_n(s) - \sigma_{n-1}(s) \bigr | < d(s)/ 2^{n+1}$ 
\item[(v)]  $\bigl\| \sigma_n |_\ab - \sigma_{n-1} |_\ab \bigr\|_{C^n} < \delta_j / 2^n$
\item[(vi)]  $X_n$ is contained in the $(1/n)$-neighborhood of $\sigma_n(\I) \cup \h E$
\item[(vii)]  $X_n$ contains the $\alpha_n$-neighborhood of $\sigma_n(\I) \cup \h E$
\end{enumerate}
}

Before constructing the sequences $(\sigma_n)$, $(X_n)$, and $(\alpha_n)$, we show that their existence will yield the theorem.  Condition (iii) implies that the sequence $(\sigma_n)$ converges uniformly to a continuous map $\sigma: \I \rightarrow \C^N$.  By condition (i), $\sigma|_K=\sigma_0|_K$, and hence $\sigma(K)=E$.  Condition (v) implies that for each $j=1, 2, \ldots$, the sequence $(\sigma_n|_\ab)$ converges in $C^k$-norm for every $k=1, 2, \ldots$.  Consequently, $\sigma|_\ab$ is of class $C^\infty$.  In addition, condition (v) implies that $\bigl\| \sigma |_\ab - \sigma_0 |_\ab \bigr\|_{C^1} < \delta_j$, and hence $\sigma |_\ab$ is a $C^\infty$-embedding.

We next verify that $\sigma$ is injective and hence is a topological embedding.  Since $\sigma |_K = \sigma_0 |_K$ is injective and as just noted each map $\sigma |_\ab$ is injective, it is enough to verify that $\sigma(s) \neq \sigma(t)$ for $s\in \openab$ and $t\in \I \sm \openab$.  Condition (iv) implies that $| \sigma(s) - \sigma_0(s) | < d(s)/2$.  If $t\in K$, then $\sigma(t)=\sigma_0(t)$ and 
$| \sigma_0(s) - \sigma_0 (t) | \geq d(s)$, so $\sigma(s) \neq \sigma(t)$.
If $t\in (a_k, b_k)$ for some $k\neq j$, then
$| \sigma_0(s) - \sigma_0(t) | \geq \max\{ d(s), d(t) \}$ and
$| \sigma(t) - \sigma_0(t) | < d(t)/2$, so

\begin{equation*}
\begin{split}
\bigl |\sigma(s) - \sigma(t) \bigr |
&\geq \bigl |\sigma_0(s) - \sigma_0(t) \bigr | - \bigl |\sigma(s) - \sigma_0(s) \bigr | - \bigl |\sigma(t) - \sigma_0(t) \bigr | \\ 
&> \max \{ d(s), d(t) \} - d(s)/2 - d(t) / 2\\
& \geq 0.
\end{split}
\end{equation*}
Thus again $\sigma(s) \neq \sigma(t)$.
Consequently, $\sigma$ is injective.

Recall that we set $J=\sigma(\I)$.
Conditions (iii) and (vii) together imply that $X_n$ contains $\sigma_m(\I) \cup \h E$ for all $m\geq n\geq0$.
The reader can verify that this and condition (vi) together imply that 
$J \cup \h E=\bigcap_{n=0}^\infty X_n$.  Since the intersection of a collection of polynomially convex sets is polynomially convex, we conclude that $J \cup \h E$ is polynomially convex.  The desired equality $\h J = J \cup \h E$ follows.
Of course the desired inclusion $J\subset \Omega$ also follows from the equality
$J \cup \h E=\bigcap_{n=0}^\infty X_n$ since $X_0\subset \Omega$.
This concludes the verification that the existence of the sequences $(\sigma_n)$, $(X_n)$, and $(\alpha_n)$ will yield the theorem.

Finally, we construct the sequences $(\sigma_n)$, $(X_n)$, and $(\alpha_n)$.  We already have $\sigma_0$, $X_0$, and $\alpha_0$.
We proceed now by induction.  Suppose for some $k\geq 0$ we have chosen
$\sigma_0,\ldots, \sigma_k$, $X_0,\ldots, X_k$, and $\alpha_0,\ldots, \alpha_k$ so that conditions (i)--(vii) are satisfied for all $n=1,\ldots, k$ (and all $j=1,2,\ldots$ and all $s\in \I \sm K$).  The set $\sigma_k^{-1}\bigl(\C^N \sm \Int(L_{k+1})\bigr)$ is a compact subset of $\I \sm K$ and hence is contained in the union $\openabn s$ for some finite $s$.  Let $d_k$ denote the minimum of the function $d$ on $\sigma_k^{-1}\bigl(\C^N \sm \Int(L_{k+1})\bigr)$.  Note that $d_k>0$.  Let $M= \abn s$.  Applying Theorem~\ref{AWold-theorem} to the compact manifold $M$, the $C^\infty$-embedding $\sigma_k|_M$, and the compact polynomially convex set $L_{k+1}$ gives that there exists a $C^\infty$-embedding $g: M\rightarrow \C^N$ such that
{\baselineskip=16pt
\begin{enumerate}
\item[(1)] $\| g - \sigma_k|_M\|_{C^{k+1}} < \min \{ \alpha_k/2^{k+1}, d_k/2^{k+2}, \delta_1/2^{k+1},\ldots, \delta_s/2^{k+1}\}$
\item[(2)] $g=\sigma_k$ on $M\cap\sigma_k^{-1}(L_{k+1})$
\item[(3)] $g(M) \cup L_{k+1}$ is polynomially convex.
\end{enumerate}
}
Define $\sigma_{k+1}:\I \rightarrow \C^N$ to coincide with $\sigma_k$ on 
$\sigma_k^{-1}(L_{k+1})$ and to coincide with $g$ on $M$.  Then $\sigma_{k+1}$ is continuous and of course $\sigma_{k+1}|_K=\sigma_k|_K=\sigma_0|_K$ and each restriction $\sigma_{k+1}|_\ab$ is of class $C^\infty$.  By condition (1), 
$\| \sigma_{k+1} - \sigma_k \|_\infty < \alpha_k / 2^{k+1}$ and 
$| \sigma_{k+1}(s) - \sigma_k(s)| < d(s) / 2^{k+2}$ for each $s\in \I \sm K$.  
Also by condition (1), $\bigl\| \sigma_{k+1} |_\ab - \sigma_k |_\ab \bigr\|_{C^{k+1}} < \delta_j / 2^{k+1}$ for each $j$.
Furthermore, since $\sigma_{k+1}(\I) \cup \h E$ is contained in $g(M) \cup L_{k+1}$, 
condition (3) gives that the polynomial hull of $\sigma_{k+1} (\I) \cup \h E$ is contained in the $1/(k+1)$-neighborhood of 
$\sigma_{k+1}(\I) \cup \h E$.
Consequently, there exists a compact polynomially convex neighborhood of $\sigma_{k+1}(\I) \cup \h E$, which we will take as $X_{k+1}$, contained in the $1/(k+1)$-neighborhood of 
$\sigma_{k+1}(\I) \cup \h E$.  The set $X_{k+1}$ contains the $\alpha$-neighborhood of $\sigma_{k+1}(\I) \cup \h E$ for some $\alpha > 0$.  Set $\alpha_{k+1} = \min \{\alpha, \alpha_k\}$.  
Now conditions (i)--(vii) hold for all $n=1,\ldots, k+1$.
This completes the induction and the proof.
\epf

\bpf[Proof of Corollary~\ref{corollary}]
Let $J$ be an arc or simple closed curve obtained from $E$ as in Theorem~\ref{general-theorem}.  That $P(J)=C(J)$ follows from \cite[Theorem~1.6.8]{Stout} (or alternatively, \cite[Lemma~3.1]{AIW2001}).
\epf

The proof of Theorem~\ref{poly-convex} will use the following 
lemma.  The special case of the lemma in which $Y$ is the closed unit ball and $x_0$ is the origin is proved in \cite{Izzo2} as \cite[Lemma~3.2]{Izzo2}.  The proof of the general case is essentially the same, but at the suggestion of a referee, we include it here for completeness.

\blem\label{0}
Let $Y$ be a compact polynomially convex set in $\C^N$, let $x_0$ be a point of $Y$, and let $\vep>0$.  Let $\{p_j\}$ be a countable collection of polynomials on $\C^N$ such that $p_j(x_0)\neq 0$ for each $j$.  Then there exists a polynomially convex Cantor set $E$ with $x_0\in E \subset Y\subset \C^N$ such that 
\item{\rm(i)} each $p_j$ is zero-free on $E$
\item{\rm(ii)} $m(Y\sm E)<\vep$.
\elem

\bpf
By multiplying each $p_j$ by a suitable complex number if necessary, we may assume that $\Re p_j(x_0)\neq 0$ for each $j$.  Furthermore, by enlarging the collection $\{p_j\}$, we may assume that $\{p_j\}$ is dense in $P(Y)$.  

For each $j$, the set $p_j^{-1}(\{z\in \C: \Re z=0\})$ is a real-analytic variety in $\C^N$ and hence has $2N$-dimensional measure zero.  Consequently, we can choose $0<\vep_j< \min\{ |\Re p_j(x_0)|, 1\}$ such that 
$$m\Bigl(p_j^{-1}\bigl(\{z\in \C: |\Re z|< \vep_j\}\bigl) \cap Y \Bigr)< \vep/2^j.$$
Then
$$m\Bigl(\bigcup_{j=1}^\infty p_j^{-1}\bigl(\{z\in \C: |\Re z|< \vep_j\}\bigl) \cap Y \Bigr)< \sum_{j=1}^\infty \vep/2^j = \vep.$$
Thus setting $X=Y \sm \bigcup_{j=1}^\infty p_j^{-1}\bigl(\{z\in \C: |\Re z|< \vep_j\}\bigl)$ we have that $m(Y \sm X)<\vep$.  Obviously $x_0\in X \subset Y$, each $p_j$ is zero-free on $X$, and $X$ is compact.  

For each $j$, choose a closed disc $\od_j$ containing $p_j(Y)$.  Then
$$X=\bigcap_{j=1}^\infty \Bigl(p_j^{-1}\bigl(\{z\in \od_j: |\Re z| \geq \vep_j\}\bigr) \cap Y\Bigr).$$
Each set $\{z\in \od_j : | \Re z| \geq \vep_j\}$ is polynomially convex since it has connected complement in the plane.  Hence each set 
$p_j^{-1}\bigl(\{z\in \od_j: |\Re z| \geq \vep_j\}\bigl) \cap Y$ is polynomially convex (by the elementary \cite[Lemma~3.1]{Izzo2018}).  Consequently, $X$ is polynomially convex.

Let $x$ and $y$ be arbitrary distinct points of $X$.  Because $\{p_j\}$ is dense in $P(Y)$, there is some $p_j$ such that $\Re p_j(x)>1$ and $\Re p_j(y)<-1$.  Because $| \Re p_j| \geq \vep_j$ everywhere on $X$, it follows that $x$ and $y$ lie in different components of $X$.  Consequently, $X$ is totally disconnected.

By the Cantor-Bendixson theorem \cite[Theorem~2A.1]{Mos} the compact metric space $X$ is the disjoint union of a perfect set $E$ and an at most countable set.  Then $m(Y \sm E)=m(Y \sm X)$, and $E$ is a Cantor set by the usual characterization of Cantor sets as the compact, totally disconnected, metrizable spaces without isolated points.  Finally, $E$ is polynomially convex, for $\h E$ is contained in $X$, so if $\h E$ were strictly larger than $E$, then $\hh E$ would have isolated points, which is impossible (by the Shilov idempotent theorem, for instance).
\epf

\begin{proof}[Proof of Theorem~\ref{poly-convex}]
We treat only the construction of the arc, the construction of the simple closed curve being essentially the same.
Because $\Omega$ is a Runge domain of holomorphy, $\Omega$ can be exhausted by compact polynomially convex sets.  Consequently, there is a compact polynomially convex set $Y$ such that $x_0\in Y \subset \Omega$ and $m(\Omega \sm Y) < \vep/2$.
Using Theorem~\ref{general-theorem} and Lemma~\ref{0}, we will obtain the desired arc by an argument similar to the proof of \cite[Lemma~3.3]{Izzo2}.
Choose a countable collection $\{p_j\}$ of polynomials that is dense in $P(\overline\Omega)$ and such that $p_j(x_0)\neq 0$ for each $j$.  By Lemma~\ref{0}, there exists a polynomially convex Cantor set $E$ such that $x_0\in E \subset Y$, each $p_j$ is zero-free on $E$, and $m(Y\sm E) < \vep/2$.  Let $J$ be an arc in $\Omega$ obtained from $E$ as in Theorem~\ref{general-theorem}.  Then of course $J$ is polynomially convex and $m(\Omega \sm J) < \vep$.  As noted in the proof of Corollary~\ref{corollary}, $P(J)=C(J)$.
Let $\sigma: [0,1]\rightarrow J$ be a homeomorphism.  Assume without loss of generality that $\sigma(0)$ and $\sigma(1)$ are in the Cantor set $E$, and let $(a_1,b_1)$, $(a_2,b_2)$, $\ldots$ be the disjoint open intervals whose union is $[0,1]\sm \sigma^{-1}(E)$.  Because each $p_j$ is zero-free on $E$, there exists 
$\vep_j>0$ such that $p_j(E)$ is disjoint from $\{z\in \C: |z|<\vep_j\}$.  Each set $p_j([a_k, b_k])$ has empty interior in the plane, so it follows that the set 
$p_j(J)\cap \{z\in \C: |z|<\vep_j\}$ is a countable union of nowhere dense set and thus has empty interior.  Thus there exist arbitrarily small complex numbers $\alpha$ such that $p_j+\alpha$ has no \hbox{zeros} on $J$.  Consequently, the set of polynomials zero-free on $J$ is dense in $P(\overline\Omega)$.
\end{proof}

\begin{proof}[Proof of Theorem~\ref{arc-no-disc}]
By \cite[Theorem~7.1]{Izzo} there exists a Cantor set $E$ in $\C^3$ such that $\h E \sm E$ is nonempty and $P(E)$ has dense invertibles.  Let $J$ be an arc or simple closed curve obtained from $E$ as in Theorem~\ref{general-theorem}.  Then $\h J \sm \h E$ has two-dimensional Hausdorff measure zero.  Consequently, by \cite[Corollary~1.6.8]{Stout} for instance, $P(\h J)=\{ f\in C(\h J): f|_{\h E}\in P(\h E)\}$.
Therefore, the density of invertible elements in $P(E)$ implies that $P(J)$ has dense invertibles by \cite[Lemma~8.1]{Izzo}.  Furthermore, $\h J\sm J=\h E\sm J$ and the set $\h E\sm J$ must be nonempty as the two-dimensional Hausdorff measure of $\h E\sm E$ can not be zero, by \cite[Corollary~1.6.8]{Stout} for instance (and in fact must be infinite by \cite[Theorem~21.9]{AW}).
\end{proof}


\begin{thebibliography}{BC87}

\bibitem{A} H. Alexander, {\it Polynomial approximation and hulls in sets
of finite linear measure in ${\bf C}^n$},
Amer.\ J.\ Math.\ {\bf 93} (1971), 65--74.

\bibitem{AW} H. Alexander and J. Wermer, {\it
Several Complex Variables and Banach Algebras}, 3rd~ed.,
Springer, New York, 1998.

\bibitem{AIW2001} J. T. Anderson, A. J. Izzo, and J. Wermer, {\it Polynomial approximation on three-dimensional real-analytic submanifolds of $\mathbb{C}^n$},  \emph{Proc.\ Amer.\ Math.\ Soc.\ }{\bf 129} (2001) no. 8, 2395--2402.

\bibitem{Antoine} L. Antoine, {\it Sur l'hom\'eomorphie de deux figures et de leurs voisinages}, J. Math.\ Pures Appl.\ (8) vol. 4 (1921), 221--325.

\bibitem{AWold} L. Arosio and E. F. Wold, {\it Totally real embeddings with prescribed polynomial hulls}, Indiana Univ.\ Math.\ J. {\bf 68\/} (2019), 629--640.

\bibitem{DalesF}
H. G. Dales and J. F. Feinstein, {\it Banach function algebras with dense invertible group\/}, Proc.\ Amer.\ Math.\ Soc., {\bf 136\/} (2008), 1295 - 1304.

\bibitem{Izzo2018} A. J. Izzo,
{\it Gleason parts and point derivations for uniform algebras with dense invertible group},
Trans.\ Amer.\ Math.\ Soc.\/ {\bf 370} (2018), 4299--4321.

\bibitem{Izzo} A. J. Izzo, {\it Spaces with polynomial hulls that contain no analytic discs\/}, Math.\ Ann.\ (accepted).

\bibitem{Izzo2} A. J. Izzo, {\it Gleason parts and point derivations for uniform algebras with dense invertible group II\/} (submitted).

\bibitem{ISW} A. J. Izzo, H. Samuelsson Kalm, and E. F. Wold, {\it Presence or absence of analytic structure in maximal ideal spaces},
Math.\ Ann.\ {\bf 366} (2016), 459--478.

\bibitem{IS} A. J. Izzo and E. L. Stout, {\it Hulls of surfaces}, Indiana Univ.\ Math.\ J. {\bf 67} (2018), 2061--2087.

\bibitem{Mos} Y. N. Moschovakis,  {\it Descriptive Set Theory\/}, Studies in Logic and the Foundations of Mathematics, {\bf 100\/}, North-Holland Publishing Co., Amsterdam-New York, 1980.

\bibitem{Rudin} W. Rudin, {\it Subalgebras of spaces of continuous functions},
Proc.\ Amer.\ Math.\ Soc.\ {\bf 7} (1956), 825--830.

\bibitem{Stol1}  G. Stolzenberg, {\it A hull with no analytic structure\/},
J.\ Math.\ Mech.\  {\bf 12} (1963), 103--111.

\bibitem{Stout}  E. L. Stout,  {\it Polynomial Convexity\/},
Birkh\"auser, Boston, 2007.

\bibitem{Wermer1954} J. Wermer,  {\it Algebras with two generators\/},
Amer.\ J. Math.\ {\bf 76} (1954), 853--859.

\bibitem{Wermer-Cantor}  J. Wermer,  {\it Polynomial approximation on an arc in $C^3$\/},
Ann of Math. (2) {\bf 62\/} (1955), 269--270.

\bibitem{Wermer1958} J. Wermer, {\it The hull of a curve in ${\bf C}^n$},
 Annals Math.\ {\bf 68} (1958), 550--561.

\bibitem{Whyburn} G. T. Whyburn, {\it Concerning the proposition that every closed, compact, and totally disconnected set of points is a subset of an arc\/}, Fundamenta Matematica {\bf 18\/} (1932), 47--60.



\end{thebibliography}
\end{document}